\documentclass[bib]{ba}

\usepackage{epsfig}
\usepackage{bayes}

\usepackage{amsmath,amsfonts,amsthm,amssymb,latexsym}
\usepackage{natbib}

\newtheorem{example}{Example}[section]
\newcounter{exA}

\newcommand\dropcap\noindent

\begin{document}

\inserttype{article}
\author{Robert, C.P.}{%
  {\sc Christian P.~Robert}\\Universit\'e Paris-Dauphine, CEREMADE, and CREST, Paris
}	
\title[Reading The Search for Certainty]{{\it{\bfseries The Search for Certainty}\/}:\\ A critical assessment}

\maketitle

\begin{abstract}
{\em The Search for Certainty} was published in \citeyear{burdzy:2009} 
by Krzysztof Burdzy. It examines the {\em ``philosophical duopoly"} of von Mises and de Finetti
at the foundation of probability and statistics and find this duopoly missing. This review exposes
the weakness of the arguments presented in the book, it questions the relevance of introducing 
a new set of probability axioms from a methodological perspective, and it concludes at the lack of 
impact of this book on statistical foundations and practice.
\end{abstract}


\section{Introduction}
    \begin{quote}{\em 
    ``This book is about one of the greatest intellectual failures of the twentieth century---several unsuccessful 
      attempts to construct a scientific theory of probability." {\em The Search for Certainty}, page vii.
    }\end{quote}
Intrigued by the premises and claims found on the back-cover of the book, I read through Krzysztof Burdzy's {\em The Search for Certainty}, 
expecting a novel analysis of the philosophical foundations of statistics. Instead, I found a highly focussed
discussion on two philosophical approaches to the definition of probability that were reminiscent of the
discussions found at the beginning of the 20th Century and I thus closed the book 
concluding that it bears little if any connection with our field, while acknowledging the original
perspective of the author. The present paper is the summary of my perusal of the book, pointing out shortcomings of the book 
from a statistician's perspective. 

The entry reproduced above (which happens to be the first sentence of the book) illustrates how high the author sets his goal, 
but {\em The Search for Certainty} does not substantiate the ultimate generality of this claim. Indeed,
as exposed more clearly by the subtitle {\em On the Clash of Science and Philosophy of Probability}, the focus of the book 
is a radical (if, in my opinion, superficial) criticism of the philosophical foundations of probability and (indirectly) of 
statistics, rather than a re-analysis of probability calculus (in the sense of \citealp{billingsley:1995}), which remains untouched 
throughout the book---even though Section \S 11.4 briefly covers Kolmogorov's axiomatic system with the valid comment that it does 
not represent a philosophical theory. In addition, {\em The Search for Certainty} does not cover the impact of such philosophical theories 
on science. As argued below, I consider that the book fails short of providing a convincing argument of how probability concepts 
should be modified to enhance the role of probability in science in general and in statistics in particular. While not a philosopher 
myself, and thus barred from making comments about the philosophical validity of the book {\em per se}, I am under the impression that 
the depth of scholarship argumentation in {\em The Search for Certainty} in terms of philosophy of sciences and even of epistemology 
is not sufficient to make a case for this perspective. Despite the author's claim that ``this book contains a 
number of ideas that I have not seen anywhere in any form" (p.viii), my own feeling
is that it does not open a new direction for handling probability in science. In 
particular, when the author states that ``the book is concerned with the substance of philosophical claims and their relation with 
statistics'' (p.10), the
only conclusion I find appealing is that ``the original theories of von Mises and de Finetti are completely unrelated to statistics" (p.12)
as their philosophical theories have not had a lasting impact on statistics. 

\section{A lack of foundations}
    \begin{quote}{\em 
    ``The philosophical essence of the science of probability is to generate a single prediction." {\em The Search for Certainty}, page 52.
    }\end{quote}
The central ``philosophical object" in the book is the notion of the probability of a single event (not to be confused with the mathematical
definition). Because this event cannot be reproduced, the author argues that this rules out an approach based on an infinite sequence of 
events, hence also ruling out the frequentist (i.e., von Mises') definition,
while a choice of a formally valid probability distribution covering this event is subjective, thus non-verifiable and hence un-scientific.
As often the case in modern debates about the nature of induction, the author calls for the scholarly authority of Karl Popper, 
``the champion of the propensity theory of probability" (p.43), and of his falsifiability criterion, using de Finetti's statement 
that ``no probability statement is falsifiable in any sense" (p.22), to support his own arguments. In my opinion,
this is the extent of the theoretical support for the criticisms contained in the book about the failure of both 
von Mises' and de Finetti's theories. {\em The Search for Certainty} contains no philosophical background other than this call to 
Popper---who, incidentally, also demonstrated in \cite{popper:miller:1983} 
the unfalsifiability of statistical induction\footnote{The author seems to follow the same track when defining
{\em predictions} as events with very high probabilities, 
a definition which obviously removes most of the problems linked with induction. This may be inspired from 
\cite{definetti:1974} who distinguishes between {\em previsions} and {\em predictions} (p.134), as well as
Popper's proposal to only consider events with probabilities close to $0$ or $1$, reproduced on page 18.} as a whole.
(This is furthermore acknowledged as a ``repackaging of Popper's idea for general consumption", p.8.)
A few additional references like Hofstadter's vulgarisation book about G\"odel's incompleteness theorem are mentioned in Chapter 3,
but the lack of connections with the existing literature on the philosophy of science and epistemology appears to me as an
indicator of the isolated position of the author within those domains. 

The tone set in the book resembles those found in probability books from the turn of the previous century, that is between the 
19th and the 20th Centuries (as well as in \citealp{keynes:1921} and \citealp{jeffreys:1931,jeffreys:1939}) about the logical
nature of probability. For instance, a common feature is that the author proposes new axioms to replace Kolmogorov's classical 
three axioms (because those ``say nothing about how to match the mathematical results with reality" (p.31)):
\begin{itemize}
\item[(L1)] Probabilities are numbers between $0$ and $1$, assigned to events whose outcomes may be unknown.
\item[(L2)] If events $A$ and $B$ cannot happen at the same time then the probability that one of them will occur is the sum of probabilities of the
individual events, that is, $P(A\text{ or }B) = P(A) + P(B)$.
\item[(L3)] If events A and B are physically independent then they are independent
in the mathematical sense, that is, $P(A \text{ and }B) = P(A)P(B)$. 
\item[(L4)] If there exists a symmetry on the space of possible outcomes which maps an event $A$ onto an event $B$ then the two events have equal probabilities, that is, $P(A) = P(B)$.
\item[(L5)] An event has probability $0$ if and only if it cannot occur. An event has probability $1$ if and only if it must occur.
\end{itemize}
It is hard both to come up with satisfying epistemological justifications for 
introducing in 2009 a new set of axioms (even though the author maintains they ``are implicit in all textbooks" (p.35))
and not to find difficulties with them. The second part of axiom (L1) is novel, but it involves an observer and thus implies endless
complications {\em (what if there are two observers? an infinity of them? how does the observer impact on the event?}, etc.)
Axiom (L2) is your standard Kolmogorov's mutually-exclusive-events second axiom. Axiom (L3) is usually a definition with no connection to 
physical properties, which here involves an observer and further hints at a basic determinism that pervades the whole book (as illustrated by
``nobody knows the objective truth", p.164). From a Bayesian perspective, axiom (L4) relates 
to both Laplace's Principle of Insufficient Reasons \citep{keynes:1921,stigler:1986}, 
whose limitations have been the cause of much debate in the selection of prior distributions, and 
to invariance principles that lead to Haar measures as default noninformative priors \citep{berger:1985}. This notion
of symmetry is first understood as invariance in the book (\S 3.1) but it later turns into a wider and less defined notion (see, e.g., the
discussion of the connection between symmetry and priors in \S 8.4.3), while the connection with Haar measures and invariant priors is never
made. In any case, this axiom once again involves an observer, since Burdzy maintains that symmetry is ``relative" (\S 3.5), as in independence
versus exchangeability assumptions (\S 3.14). At last, (L5) is formulated in a psychological wording rather than in the mathematically 
unambiguous constraints $P(\emptyset)=0$ and $P(\Omega)=1$.

Given that the author bases his whole argumentation on this set of axioms, his examples are almost necessarily limited to the coin-and-urn 
types, ruling out more complex models such as those found in continuous spaces. Burdzy concedes as much when stating that ``there is no presumption 
that it should be trivial to derive popular models such as linear regression or geometric Brownian motion" (p.38). I would have been more
open to the relevance of this change of axioms had involved examples (for modelling or inference) been analysed within the book. The absence
of such examples seems to signal a fundamental difficulty in connecting those philosophical reflections with our field.

\section{Which relevance for statistics?}
\subsection{Classical statistics}
    \begin{quote}{\em ``Mr. Winston is unique because we know something about him that we do not know about 
    any other individual in the population." {\em The Search for Certainty}, page 66.}\end{quote}
A lack of convincing arguments transpires from the study of statistics found in the book, starting with the 
above quote that simply expresses the fact that if we condition on too many covariates there is no predictive ability left 
for the model. The criticisms on decision theory found in Chapter 4 are in my opinion shallow, ranging from 
questioning the purpose of maximising the expected gain (``why not minimize the third moment of the gain?", p.78; 
``maximizing the subjective gain is tautological", p.79) 
to the fact that ``not all decisions are comparable" (p.82), 
to the application of utility theory to non-monetary rewards (\S 4.4.3),
to the issue that utility is non-linear for monetary rewards (\S 4.4.2), but with no mention made of much earlier studies 
of this problem like Laplace's Saint-Petersburg paradox \citep{laplace:1795}, 
to the impossibility of aggregated actions, again with no mention made of Allais' 
(\citeyear{allais:1953}) paradox and of the subsequent literature,
to the paradoxical objection that it operates ``in a
situation involving uncertainty when no relevant information is available" (p.86) and should thus offer ``no scientific advice" (p.87).

Chapter 6 on classical statistics is also very reductive in that confidence intervals are processed
in the sense of Burdzy's predictions, which has a ``long list of practical problems" (p.118), 
in association with the insistence that predictions have to be verified---thus contradicting
both the initial perspective of single occurrence events and the purpose of statistics. We find there the unusual definition that the ``goal 
of the estimation theory is to find an explicit formula for the distribution of a given estimator" (p.119) and the assimilation between 
frequentist estimation and unbiasedness, missing the important feature that unbiased estimators only exist in a very limited number of 
cases \citep{lehmann:casella:1998}, and leading to the weak argument that unbiasedness is not supported by long term frequency properties. 
Surprisingly, (classical) hypothesis testing is getting good grades: ``the theory of hypothesis testing is a well-justified science" 
(p.124), even though true predictions in Burdzy's sense would imply almost always accepting the null hypothesis since they push the 
type I error all the way to zero. The only criticism of classical hypothesis testing advanced by Burdzy
is paradoxically that it ``involves not only the `null 
hypothesis' but also an alternative hypothesis" (p.129) which is a standard argument used against Bayesian testing 
\citep{templeton:2008,clade:2010}. Furthermore, no mention is made of maximum likelihood estimation in this range of 
criticisms, Fisher being only perceived as an opponent to Jeffreys (p.11 and p.68)---who, incidentally, is presented 
as a ``subjectivist" (p.245).

\subsection{Bayesian statistics}
    \begin{quote}{\em ``Bayesian statisticians see nothing wrong with collecting data first and starting the statistical 
    analysis later." {\em The Search for Certainty}, page 179.}\end{quote}
The bases of Bayesian statistics are also misrepresented in that, while the author's position is overall sympathetic (``Bayesian statistics is
a very successful branch of science because it is capable of making excellent predictions", p.177), DeGroot's (1970), as well as Savage's (1972), 
axiomatic derivations of a prior distribution from a consistent set of axioms are criticised on the basis that they do ``not guarantee a success in 
practical life" (p.25). Most of Chapter 8 on Bayesian statistics operates at the level of the ambiguous meaning of the word {\em subjective}, a
debate I do not find of particular relevance. For instance, the attempt at distinguishing between de Finetti's theory of priors,
which are ``a complete probabilistic representation of the universe" (p.179),
and Bayesian model-based priors is unclear, since the next paragraph makes mentions de Finetti's representation of exchangeable sequence 
as i.i.d.~sequences over a mixing distribution, which obviously works as an implicit prior. Asserting that subjectivists could maintain that 
``no matter what happens, nothing will prove that any
particular model is wrong" (p.180) is reductive, while the insistence on priors being ``chosen after collecting the data" (p.181) is 
misrepresenting the Bayesian practice. Another misconception stands with objective priors which should ``involve probability assignments
that can be objectively verified" (p.182), as if there was a true prior at the ready.
The alternative vision of Bayesian statistics as an iterative method (\S 8.4.2) is simply asserting the coherence of the posterior 
update within a probabilistic setting (and has no connection with the {\em prior feedback} technique of \citealp{robert:1993}) and 
this analysis misses the overall consistency properties of most Bayesian procedures \citep{ibragimov:hasminskii:1981,diaconis:freedman:1986}. 
In fact, Burdzy seems to reject asymptotic consistency as unrealistic and he argues 
that, for a given sample size, there always are two priors that differ enough for the posteriors to differ as well. While technically
true, this cannot be seen as a deterrent against consistency.

Insisting on ``the ultimate criterion for the choice of a prior [being] the
reliability of predictions generated by the posterior distribution" (p.185) shows a deep misunderstanding of 
statistics. (A more generous view would be to relate this perspective to the choice of reference priors, \`a la 
\cite{bernardo:1979}, and of matching priors, as in \cite{welch:peers:1963}, but neither is mentioned in the book.) 
Burdzy goes on questioning DeGroot's (\citeyear{degroot:1970}) choice of title as ``the optimality of your decisions 
is tautological" (p.187). This is missing the whole array of theoretical works started by Abraham Wald in the 1950's 
about the admissibility of Bayesian procedures, complete class theorems and the like 
\citep[see, e.g.,][]{berger:1985,robert:2001}, as well as the notion of minimal regret well-covered by \cite{savage:1954}. 
Playing devil's advocate, I further find surprising that, from a philosophical point of view, the author accepts so readily the choice
of priors by Bayesian statisticians and dismisses their partial arbitrariness, given that this is usually the 
broadest entry point for criticisms of the Bayesian approach \citep{gelman:2008}. Similarly, the book does not 
mention at all the delicate calibration of Bayes factors \citep{jeffreys:1939} or the similarly delicate
handling of improper priors, which are also morsels of choice for Bayesian critics \citep{robert:chopin:rousseau:2009}.

The book further contains a (short) philosophical assessment of \cite{berger:1985} 
and \cite{gelman:carlin:stern:rubin:2001} in Section 8.10, blaming the latter for mentioning subjectivity at all, 
and the former for using ``too many philosophical arguments" (p.128), because one justification (rather than
seven, \citealp{berger:1985}, \S 4.1) should be ``sufficient" (p.194). This does not constitute a proper level 
of academic philosophical argumentation. In discussing the issue of ``competing non-informative priors" (p.193), which is genuine, Burdzy 
once again calls for ``verifying predictions" (p.194), misunderstanding the purpose of statistics. Nonetheless, he
agrees on the point that the ``distinction [between constants and random variables is] irrelevant" (p.195), taking 
as an example the speed of light already treated in \cite{berger:1985}.

\section{Narrow focus}
    \begin{quote}{\em ``Of the four well crystallised philosophies of probability, two chose the certainty as their intellectual holy grail. 
    Those are the failed theories of von Mises and de Finetti." {\em The Search for Certainty}, page 30.}\end{quote}
What the author perceives as the main point in {\em The Search for Certainty}, at least according both to the preface and to the book back-cover, 
reproduced in the quote above,  is very narrow and thus of restricted interest for Bayesian statistics.
Indeed, Burdzy concentrates on two very specific entries to frequentism and subjectivism, namely von Mises' and de Finetti's, respectively, 
while those are not your average statistician's references. 
For instance, \cite{vonmises:1957} bases his definition of frequency properties on the notion of {\em collectives}, a notion I had not 
previously encountered. The criticisms in Chapter 5 are not more compelling than those found in other chapters, see for 
instance the point that collectives cannot be found in reality. (An interesting connection is made, though, with 
Marsaglia's {\em die hard} battery of tests used in checking random generators  see \citealp{robert:casella:2004}.) Or the other point that the 
Law of Large Numbers is ``impossible to apply in real life" (p.107).

    \begin{quote}{\em ``There is no justification for the use of the Bayes theorem in the subjective theory." 
    {\em The Search for Certainty}, page 144.}\end{quote}
Similarly, de Finetti's statement ``Probability does not exist" is over-exploited in the book. Obviously, this cannot be 
seen as the core principle for most Bayesian statisticians and few still relate to his all-subjective 
stance for conducting Bayesian inference (a fact acknowledged in Chapter 8). Chapter 7 on the analysis of de Finetti's theory is therefore
relevant only for a very limited range of researchers and even those interested in the topic for epistemological or 
philosophical reasons may find the corresponding discussion too superficial to be informative. The tone of this chapter is 
in addition rather forceful, like stating that ``de Finetti's theory fails to report any probabilistic facts
observed in the past" (p.137), which seems to contradict the updating principle represented by the Bayes theorem and
analysed previously in \S 8.4.2. While I (and \citealp{berger:1985}, among others) agree that ``consistency alone is 
[not] a sufficient basis" (p.139) for supporting a procedure or a theory, I fear the simplified representation 
of de Finetti's approach produced in {\em The Search for Certainty} leads to rather
surprising statements like the one above or like the point that updating the prior based on past observations is useless 
(``the subjective postulates implicitly tell the statistician that they can completely ignore the data", p.149). Having 
reached this extreme position makes it very difficult to argue any further, and I will thus not cover the remaining arguments 
in the chapter.  
(The chapter repeatedly mentions problems linked to what appears to be a change of
model along observations, as in \S 7.6.4 and \S 7.6.6, but this is not clear enough to argue, in my opinion.) The only part I
found myself agreeing with was Section \S 7.15 on arbitrage, since Black-Scholes option pricing
is indeed a mathematical construct with no connection with reality.

\section{A matter of style}
    \begin{quote}{\em ``The expected value is hardly ever expected." {\em The Search for Certainty}, page 207.}\end{quote}
While this is a minor issue, compared with the above, I find the style of {\em The Search for Certainty} to be somehow annoying
at times. First, it is predominantly non-technical, the worked-out examples always being of the coin and balls-in-an-urn type, with no
statistics example. Second, retorical arguments are never more than one paragraph long. Metaphors and weak analogies abound but 
are rarely helpful, as illustrated by the chapter on confusing terms ({\em Abuse of Language}, pp.207--209) which contains aphorisms 
like the one above. {\em Bis repetita placent}, the lack of support from other scholarly sources is somehow unsettling.

The final chapters about {\em Teaching Probability} (Chap.9), {\em What is Science?} (Chap.11), and 
{\em What is Philosophy?} (Chap.12) appear as posterior add-ons that do not contribute much to the central
debate but mostly repeat earlier arguments. (The remaining four chapters are to be considered as appendices of sort.)
The author argues that teaching the foundations of probability is in a ``confused state" while I think those are 
foundations are simply not taught at all. Apart from Italy, I do not know of any probability nor statistics course that includes either
von Mises' or de Finetti's theories in their curriculum. Exposing the students to Kolmogorov's axiomatic system appears
to me like the natural entry to a probability (calculus) course, while the author states that ``this does more harm than
good" (p.200), a controversial position about a purely mathematical concept. 

In an even narrower perspective, note that the cover of {\em The Search for Certainty} is a puzzle in itself. 
The event it produces is associated with the ``wrong" probability, involuntarily or not. Indeed, this cover shows seven dies 
with sixes on top, along with $p=0.0000036...$, which is equal to $6^{-7}$ and is indeed the probability to obtain seven times 
the same value with fair dies. However, the
dies are shown in such a way that they all display fours or/and fives on their visible sides, which means the probability 
should be something between $(6\times 2)^{-7}=2\,10^{-8}$ and  $(6\times 4)^{-7}=2\,10^{-10}$ instead. 
The signal is somehow misleading in that this picture corresponds either to a deterministic event---the dies 
were arranged this way---or to an extremely unlikely but reproducible event. 

\section{Conclusion}
    \begin{quote}{\em ``The incredibly high standards for the measurement of probability set by von Mises and de Finetti have
    no parallel in science." {\em The Search for Certainty}, page 232.}\end{quote}
In conclusion, I have found reading {\em The Search for Certainty} a mostly unrewarding experience, even though it opens vistas about
the very notion of probability. As such, I do not consider that this book makes a significant contribution 
to the foundations of statistical inference in general and of Bayesian analysis in particular. The statement that ``the above
philosophical claims are incomprehensible to all statisticians" (p.9) could extend to {\em The Search for Certainty} as well.
Once the book read, the reader is certainly left with no additional insight on how to conduct inference in a proper manner. The
exposition of the new axiomatic system (L1)--(L5) does not produce a working principle and we thus left ``miles away from 
implementing the idea" (p.232) with a vague hope that ``perhaps one day somebody will invent a philosophical theory
representing that scheme of thinking" (p.232). Having one's discipline examined by someone from another field, like a probabilist 
or a philosopher, certainly has a strong appeal and a potential for revealing new truths or questioning existing principles, but 
only if this is done at a deep enough level.

\section*{Acknowledgements}
The author's research is partly supported by the Agence Nationale de la Recherche (ANR, 212,
rue de Bercy 75012 Paris) through the 2007--2010 grant ANR-07-BLAN-0237 ``SPBayes". Communications
with Professor Burdzy helped in clarifying some factual inconsistencies in an earlier version and
comments from the Associate Editor were helpful for the presentation of this review.

\renewcommand{\bibsection}{\section*{References}}

\bibliographystyle{ba}

\end{document}